\documentclass[11pt]{article}
\usepackage{amsmath,epsf,amsthm,amssymb,amscd,mathrsfs,hyperref}

\newtheorem{theorem}{Theorem}[section]
\newtheorem{corollary}[theorem]{Corollary}
\newtheorem{lemma}[theorem]{Lemma}
\newtheorem{proposition}[theorem]{Proposition}
\newtheorem{example}[theorem]{Example}
\newtheorem{definition}[theorem]{Definition}

\newcommand{\zp}{{{\mathbb Z}_p}}
\newcommand{\qp}{{{\mathbb Q}_p}}

\begin{document}

\title{Combinatorial Congruences and $\psi$-Operators}

\author{Daqing  Wan\thanks{Partially supported by NSF. The author
thanks Z.W. Sun for helpful discussions on combinatorial
congruences}
\\ dwan@math.uci.edu \\ Department of Mathematics\\
University of California,  Irvine \\ CA 92697-3875\\
%Institute of Mathematics\\ Chinese Academy of Sciences,
%Beijing,P.R. China
}

\maketitle

\begin{abstract}
The $\psi$-operator for $(\varphi, \Gamma)$-modules plays an
important role in the study of Iwasawa theory via Fontaine's big
rings. In this note, we prove several sharp estimates for the
$\psi$-operator in the cyclotomic case. These estimates
immediately imply a number of sharp $p$-adic combinatorial
congruences, one of which extends the classical congruences of
Fleck (1913) and Weisman (1977).

\end{abstract}

\section{Combinatorial Congruences}
Let $p$ be a prime, $n\in\mathbb{Z}_{>0}$. Throughout this paper,
let $[x]$ denote the integer part of $x$ if $x\geq 0$ and $[x]=0$
if $x<0$. In the author's course lectures \cite{Wa} on Fontaine's
theory and $p$-adic L-functions given at UC Irvine (spring 2005)
and at the Morningside Center of Mathematics (summer 2005), the
following two congruences were discovered.

\begin{theorem}For integers $r\in\mathbb{Z}$,
$j\geq 0$, we have
    $$\sum_{k\equiv r({\rm mod} p)}(-1)^{n-k}\binom{n}{k}\binom{\frac{k-r}{p}}{j}\equiv 0~({\rm mod} p^{[\frac{n-1-jp}{p-1}]}).$$

\end{theorem}

We shall see that the theorem comes from a simple estimate of
$\psi(\pi^n)$ for the cyclotomic $\varphi$-module.

\begin{theorem}For integer $j\geq 0$, we have
$$\sum_{\substack{i_0+\cdots+i_{p-1}=n \\ i_1+2i_2+\cdots
\equiv r({\rm mod} p)}}\binom{n}{i_0 i_1 \cdots
i_{p-1}}\binom{\frac{i_1+2i_2+\cdots-r}{p}}{j}\equiv 0 ~({\rm mod}
p^{[\frac{n(p-1)-jp-1}{p-1}]}).$$\end{theorem}

As we shall see, this theorem comes from a simple estimate of
$\psi(\pi^{-n})$ for the cyclotomic $\varphi$-module. Note that
when $p=2$, Theorem 1.2 is equivalent to Theorem 1.1.

The above two congruences can be extended from $p$ to $q=p^a$,
where $a$ is a positive integer. To do so, it suffices to estimate
the $a$-th iterate $\psi^a(\pi^n)$. This can be done by induction.
The estimate of $\psi^a(\pi^n)$ for $n>0$ leads to

\begin{theorem} For integers $r\in\mathbb{Z}$, $j\geq 0$ and
$a>0$, we have
$$\sum_{k\equiv r({\rm mod}p^a)}(-1)^{n-k}\binom{n}{k}\binom{{k-r\over p^a}}{j}\equiv 0({\rm mod}
p^{[{n-p^{a-1}-jp^a \over p^{a-1}(p-1)}]}).$$
\end{theorem}

The estimate of $\psi^a(\pi^n)$ for $n<0$ leads to

\begin{theorem} Let
$$S_j(n,r,p^a)=\sum_{\begin{subarray}{c}i_0+\cdots+i_{p^a-1}=n\\
i_1+2i_2+\cdots \equiv r({\rm
mod}p^a)\end{subarray}}\binom{n}{i_0\cdots
i_{p^a-1}}\binom{(i_1+2i_2+\cdots -r)/p^a}{j}.$$ Then for integer
$j\geq 0$, we have
$$S_j(n,r,p^a)\equiv 0({\rm mod} p^{[{(an-a+1)(p-1)-j(ap-a+1)-1\over p-1}]}).$$
\end{theorem}

As Z.W. Sun informed me, the special case $j=0$ of Theorem 1.1.1
was first proved by Fleck \cite{Di} in 1913, and the special case
of Theorem 1.1.3 for $j=0$  was first proved by Weisman \cite{We}
in 1977. A different extension of Theorem 1.1.1 and Weisman's
congruence has been obtained by Z.W. Sun \cite{Su} using different
combinatorial arguments. Motivated by applications in algebraic
topology, Sun-Davis \cite{SD} proved yet another extension:
$$\sum_{k\equiv r({\rm mod}p^a)}(-1)^{n-k}\binom{n}{k}\binom{{k-r\over p^a}}{j}\equiv 0({\rm mod}
p^{\left( {\rm ord}_p([n/p^{a-1}]!) - j -{\rm
ord}_p(j!)\right)}).$$

\section{The operator $\psi$}

Let $p$ be a fixed prime. Let $\pi$ be a formal variable. Let $$
A^{+} =\zp [[\pi]]$$ be the formal power series ring over the ring
of $p$-adic integers. Let $A$ be the $p$-adic completion of
$A^{+}[{1\over \pi}]$, and let $B=A[{1\over p}]$ be the fraction
field of $A$. The rings $A^{+}$, $A$ and $B$ correspond to
$A^{+}_{\qp}$, $A_{\qp}$ and $B_{\qp}$ in Fontaine's theory.

We shall not discuss the Galois action on $A$, which is not needed
for our present purpose. The Frobenius map $\varphi$ acts on the
above rings by
$$\varphi(\pi) =(1+\pi)^p-1.$$
If we let $[\varepsilon] =1+\pi$, then $\varphi([\varepsilon])
=[\epsilon]^p$. The map $\varphi$ is injective of degree $p$. This
gives

\begin{proposition}
$\{1,\pi,\cdots,\pi^{p-1}\}$ (and
 $\{1,[\varepsilon],\cdots,[\varepsilon]^{p-1}\}$) is a basis of
$A$ over the subring $\varphi(A)$.
\end{proposition}

\begin{definition}
\rm The operator $\psi:A\to A$ is defined by
$$\psi(x) = \psi\left(\sum_{i=0}^{p-1}[\varepsilon]^i\varphi(x_i)\right)=x_0 =
\frac{1}{p} \varphi^{-1}({\rm Tr}_{A/\varphi(A)}(x)), $$ where $x:
A \rightarrow A$ denotes the multiplication by $x$ as
$\varphi(A)$-linear map.
\end{definition}

\begin{example}
$$\psi([\varepsilon]^n)=\left\{\begin{array}{ll}
[\varepsilon]^{n/p}, &\mbox{if }p\mid n;\\
0, &\mbox{if }p\nmid n.
\end{array}\right.$$
\end{example}
It is clear that  $\psi$ is $\ \varphi^{-1}$-linear:
$$\psi(\varphi(a)x) = a \psi(x) \ \ \forall \ a, x \in A.$$
\begin{example}
Let $a$ be a positive integer relatively prime to $p$. Then
$$\psi({1\over (1+\pi)^a-1} ) = {1\over (1+\pi)^a-1}.$$
In fact,
\begin{eqnarray*}
\psi\left(\frac{1}{[\varepsilon]^a-1}\right)&=&
\psi\left(\frac{1}{[\varepsilon]^{ap}-1}\cdot\frac{[\varepsilon]^{ap}-1}{[\varepsilon]-1}\right)\\
&=&\frac{1}{[\varepsilon]^a-1}\psi\left(1+[\varepsilon]^a+\cdots+[\varepsilon]^{(p-1)a}\right)\\
&=&\frac{1}{[\varepsilon]^a-1}=\frac{1}{(1+\pi)^a-1}.
\end{eqnarray*}
\end{example}

By $p$-adic continuity, the above example holds for any $p$-adic
unit $a\in \zp^*$. In the general theory of $(\varphi,
\Gamma)$-modules, it is important to find the fix points of $\psi$
for applications to $p$-adic L-functions and Iwasawa theory. In
the simplest cyclotomic case, we have the following description
for the fixed points (see \cite{Wa}).

\begin{proposition}
\begin{equation*}
    A^{\psi=1}=\frac{1}{\pi}\zp\oplus\zp\oplus\left\{\sum_{k=0}^{\infty}\varphi^k(x)\
    \Big | \ x\in
    \bigoplus_{i=1}^{p-1}[\varepsilon]^i \varphi(a_i + \pi
    \mathbb{Z}_p[[\pi]]), \sum_{i=1}^{p-1}a_i=0 \right\},
\end{equation*}
where  $a_i\in \mathbb{Z}_p$.
\end{proposition}

For example, if $a$ is a positive integer relatively prime to $p$,
then the element
$$\frac{a}{(1+\pi)^a-1} - \frac{1}{\pi} \in (A^+)^{\psi=1}$$
gives the cyclotomic units and the Euler system. This element is
the Amice transform of a $p$-adic measure which produces the
$p$-adic zeta function of $\mathbb{Q}$. This type of connections
is conjectured to be a general phenomenon for $(\varphi,
\Gamma)$-modules coming from global $p$-adic Galois
representations.

\section{Sharp estimates for $\psi$}

The ring $A$ is a topological ring with respect to the
$(p,\pi)$-topology. A basis of neighborhoods of $0$ is the sets
$p^kA + \pi^nA^{+}$, where $k\in \mathbb{Z}$ and $n\in
\mathbb{N}$. The operator $\psi$ is uniformly continuous. This
continuity will give rise to combinatorial congruences.

For $s\in A^{+}$, one checks that

\begin{eqnarray*}
\psi(\pi^ps)&=&\psi(([\varepsilon]-1)^ps)\\
&=&\psi(([\varepsilon]^p-1)s+pss_1)\\
&=& \pi\psi(s)+p\psi(ss_1)\in (p,\pi)\psi(sA^+).
\end{eqnarray*}
In particular,
$$\psi(\pi^pA^+)\subset
(p,\pi)A^+.$$ Thus, by iteration, we get
\begin{proposition}[Weak Estimate] Let $n\ge 0$. Then $$ \psi(\pi^nA^+)\subset
(p,\pi)^{[n/p]}A^+ = \sum_{j=0}^{[n/p]}\pi^jp^{[n/p]-j}A^+.$$
\end{proposition}

Since the exponent $[(n-jp)/p]$ is decreasing in $j$, this
proposition implies that for $x\in \pi^nA^+$, we have
$$\psi(x)=\sum_{j=0}^{\infty} a_j \pi^j, \ a_j\in \zp, \ {\rm
ord}_p(a_j) \geq [(n-jp)/p].$$
This already gives a non-trivial
combinatorial congruence. Let $r$ be an integer. Let us calculate
$\psi(\pi^n[\varepsilon]^{-r})$ in a different way.

\begin{lemma}
\rm $$\psi(\pi^n[\varepsilon]^{-r})=\sum_{j\ge
0}\pi^j\sum_{k\equiv r({\rm
mod}p)}(-1)^{n-k}\binom{n}{k}\binom{(k-r)/p}{j}.$$
\end{lemma}

\begin{proof}
Since $\pi= [\varepsilon]-1$ and $[\varepsilon] = 1+\pi$, we have
\begin{eqnarray*}
\psi(\pi^n[\varepsilon]^{-r}) &=& \psi(([\varepsilon]-1)^n[\varepsilon]^{-r})\\
&=&\psi\left(\sum_{k=0}^n(-1)^{n-k}\binom{n}{k}[\varepsilon]^{k-r}\right)\\
&=&\sum_{k\equiv r({\rm mod} p)}(-1)^{n-k}\binom{n}{k}[\varepsilon]^{(k-r)/p}\\
&=&\sum_{k\equiv r({\rm mod} p)}(-1)^{n-k}\binom{n}{k}\sum_{j\ge
0}\binom{(k-r)/p}{j}\pi^j\\
&=&\sum_{j\ge0}\pi^j\sum_{k\equiv r({\rm mod}
p)}(-1)^{n-k}\binom{n}{k}\binom{(k-r)/p}{j}.
\end{eqnarray*}
\end{proof}

Comparing the coefficients of $\pi^j$ in this equation and the
weak estimate, we get
\begin{corollary}[Weak Congruence]  Let $n\ge 0$.  We have
$$\sum_{k\equiv r({\rm
mod}p)}(-1)^{n-k}\binom{n}{k}\binom{(k-r)/p}{j}\equiv 0({\rm
mod}p^{[(n-jp)/p]}).$$
\end{corollary}

The above simple estimate is crude and certainly not optimal since
we ignored a factor of $\pi$. We now improve on it.

\begin{theorem}[Sharp Estimate I] For $n\ge 0,$ we have
$$\psi(\pi^nA^+)\subset\sum_{j=0}^{[n/p]}\pi^jp^{[{n-1-jp\over p-1}]}A^+.$$
\end{theorem}
\begin{proof}
We prove the theorem by induction. The theorem is trivial if $n\le
p-1$. Write
$$\varphi(\pi) =(1+\pi)^p-1 = \pi^p -p\pi s_1, \ s_1\in
A^{+}.$$Then,
$$\psi(\pi^p s) =\psi( (\varphi(\pi)+p\pi s_1)s) =\pi\psi(s) +p
\psi(\pi s_1s).$$ This proves that the theorem is true for $n=p$.
Let $n>p$. Assume the theorem holds for $\le n-1$. It follows that
$$\psi(\pi^nA^+)=\psi(\pi^p\pi^{n-p}A^+)\subseteq\pi\psi(\pi^{n-p}A^+)+p\psi(\pi^{n+1-p}A^+).$$
By the induction hypothesis, the right side is contained in
\begin{eqnarray*}
& &\pi\sum_{j=0}^{[(n-p)/p]}\pi^jp^{[{n-p-1-jp\over p-1}]}A^+ +
p\sum_{j=0}^{[(n+1-p)/p]}\pi^jp^{[{n-p-jp\over p-1}]}A^+\\
&=&\sum_{j=1}^{[n/p]}\pi^jp^{[{n-1-jp \over p-1}]}A^+ +
\sum_{j=0}^{[(n+1-p)/p]}\pi^jp^{[{n-1-jp\over
p-1}]}A^+.\end{eqnarray*}
\end{proof}

The function $[(n-1-jp)/(p-1)]$ is decreasing in $j$ and vanishes
for $j\geq [n/p]$. Comparing the coefficients of $\pi^j$ in the
lemma and the above sharp estimate, we deduce
\begin{corollary}[Sharp Congruence I] Let $r\in\mathbb{Z}$.
$$\sum_{k\equiv r({\rm mod}p)}(-1)^{n-k}\binom{n}{k}\binom{(k-r)/p}{j}\equiv 0({\rm mod}p^{[{n-1-jp\over p-1}]}),$$
where $j\geq 0$ is a non-negative integer.
 \end{corollary}

\begin{theorem}[Sharp Estimate II] For $n> 0,$ we have
$$\psi\left(\frac{1}{\pi^n}A^+\right)\subseteq\sum_{j=0}^{[n(p-1)/p]}\frac{1}{\pi^{n-j}}
p^{[{n(p-1)-jp-1\over p-1}]}A^+.$$
\end{theorem}
\begin{proof}Note that
$$\varphi(\pi)/\pi=\pi^{p-1}+\binom{p}{1}\pi^{p-2}+\cdots+\binom{p}{p-1}\in (\pi^{p-1},p),$$
so $(\varphi(\pi)/\pi)^n\in (\pi^{p-1},p)^n.$ Then
\begin{eqnarray*}\psi\left(\frac{1}{\pi^n}A^+\right)&=&
\psi\left(\frac{1}{\varphi(\pi)^n}\left(\frac{\varphi(\pi)}{\pi}\right)^nA^+\right)\\
&=&\frac{1}{\pi^n}\psi\left(\left(\frac{\varphi(\pi)}{\pi}\right)^nA^+\right)\\
&\subseteq&\frac{1}{\pi^n}\sum_{i=0}^np^{n-i}\psi(\pi^{i(p-1)}A^+).
\end{eqnarray*}

By Sharp Estimate I, we have
$$\psi(\pi^{i(p-1)}A^+)\subseteq\sum_{j=0}^{[i(p-1)/p]}\pi^jp^{[{i(p-1)-1-jp\over p-1}]}A^+.$$
Then,
\begin{eqnarray*}
\psi\left(\frac{1}{\pi^n}A^+\right)
&\subseteq&\sum_{j=0}^{[n(p-1)/p]}\frac{1}{\pi^{n-j}}\sum_{[jp/(p-1)]\leq i\leq n}p^{n-i+ [{i(p-1)-jp-1 \over p-1}]}A^+\\
&\subseteq&\sum_{j=0}^{[n(p-1)/p]}\frac{1}{\pi^{n-j}}p^{[{n(p-1)-jp-1
\over p-1}]}A^+.
\end{eqnarray*}
\end{proof}

\begin{corollary}[Sharp Congruence II] Let
$$S_j(n,r,p)=\sum_{\begin{subarray}{c}i_0+\cdots+i_{p-1}=n\\
i_1+2i_2+\cdots \equiv r({\rm
mod}p)\end{subarray}}\binom{n}{i_0\cdots
i_{p-1}}\binom{(i_1+2i_2+\cdots -r)/p}{j}.$$ Then integer $ j\geq
0$, we have
$$S_j(n,r,p)\equiv 0({\rm mod} p^{[{n(p-1)-1-jp\over p-1}]}).$$
\end{corollary}

\begin{proof}
\begin{eqnarray*}
& &\psi\left(\frac{1}{\pi^n}[\varepsilon]^{-r}\right)\\&
=&\frac{1}{\pi^n}\psi\left(\left(\frac{[\varepsilon]^p-1}{[\varepsilon]-1}\right)^n [\varepsilon]^{-r}\right)\\
&=&\frac{1}{\pi^n}\psi((1+[\varepsilon]+\cdots+[\varepsilon]^{p-1})^n\cdot[\varepsilon]^{-r})\\
&=&\frac{1}{\pi^n}\sum_{\begin{subarray}{c}i_0+\cdots+i_{p-1}=n\\
i_1+2i_2+\cdots \equiv r({\rm
mod}p)\end{subarray}}[\varepsilon]^{(i_1+2i_2+\cdots
-r)/p}\binom{n}{i_0\cdots
i_{p-1}}\\
&=& \frac{1}{\pi^n}\sum_{\begin{subarray}{c}i_0+\cdots+i_{p-1}=n\\
i_1+2i_2+\cdots \equiv r({\rm mod}p)\end{subarray}}\sum_{j\ge
0}\pi^j\binom{n}{i_0\cdots i_{p-1}}\binom{(i_1+2i_2+\cdots
-r)/p}{j}\\
&=& \sum_{j=0}^{\infty} {1\over \pi^{n-j}} S_j(n,r,p).
\end{eqnarray*}
The function $[(n(p-1)-jp-1)/(p-1)]$ is decreasing in $j$ and
vanishes for $j\geq [n(p-1)/p]$. Comparing the coefficients of
$1\over \pi^{n-j}$, the congruence follows.
\end{proof}

\section{Sharp estimates for $\psi^a$}

Let $a$ be a positive integer. In this section, we extend the
sharp estimates for $\psi$ to $\psi^a$.

\begin{theorem}[Sharp Estimate I] For $n\ge 0,$ we have
$$\psi^a(\pi^nA^+)\subseteq\sum_{j=0}^{[{n/p^a}]}\pi^jp^{[{n-p^{a-1}-jp^a \over p^{a-1}(p-1)}]}A^+.$$
\end{theorem}
\begin{proof}
We prove the theorem by induction on $a$. The theorem is true if
$a=1$. Assume now $a\geq 2$ and assume that the theorem holds for
$a-1$. Then,
\begin{eqnarray*}
\psi^a(\pi^n A^+)&=&\psi(\psi^{a-1}\pi^{n}A^+)\\
&\subseteq & \psi( \sum_{i=0}^{[n/p^{a-1}]}\pi^i
p^{[{n-p^{a-2}-ip^{a-1} \over p^{a-2}(p-1)}]}A^+)\\
&\subseteq & \sum_{i=0}^{[n/p^{a-1}]}\sum_{j=0}^{[i/p]}
\pi^j p^{[{n-p^{a-2}-ip^{a-1}\over p^{a-2}(p-1)}]+[{i-1-jp\over p-1}]}A^+\\
&\subseteq&\sum_{j=0}^{[{n/p^a}]}\pi^j\sum_{pj\leq i\leq
[n/p^{a-1}]}p^{[{n-p^{a-2}-ip^{a-1}\over
p^{a-2}(p-1)}]+[{i-1-jp\over p-1}]}A^+.
\end{eqnarray*}
The exponent of $p$ for a fixed $j$ is decreasing in $i$ and hence
the minimun exponent of $p$ is attained when $i=[n/p^{a-1}]$. The
minimun exponent is computed to be
$$[{n-p^{a-2}-[n/p^{a-1}]p^{a-1} \over p^{a-1}-p^{a-2}}]+[{[n/p^{a-1}]-1-jp \over p-1}] = [{n-p^{a-1}-jp^a \over
p^{a-1}(p-1)}].$$
\end{proof}

The proof of the lemma gives more general

\begin{lemma}
\rm $$\psi^a(\pi^n[\varepsilon]^{-r})=\sum_{j\ge
0}\pi^j\sum_{k\equiv r({\rm
mod}p^a)}(-1)^{n-k}\binom{n}{k}\binom{(k-r)/p^a}{j}.$$
\end{lemma}

Comparing the coefficients of $\pi^j$ in the lemma and the  sharp
estimate for $\psi^a$, we get
\begin{corollary}[Sharp Congruence I] Let $r\in\mathbb{Z}$.  Then
$$\sum_{k\equiv r({\rm mod}p^a)}(-1)^{n-k}\binom{n}{k}\binom{(k-r)/p^a}{j}\equiv 0({\rm mod}
p^{[{n-p^{a-1}-jp^a \over p^{a-1}(p-1)}]}),$$ where $j\geq 0$ is a
non-negative integer.
\end{corollary}

\begin{theorem}[Sharp Estimate II] For $n> 0$ and $a>0$,  we have
$$\psi^a\left(\frac{1}{\pi^n}A^+\right)\subseteq\sum_{j=0}^{[{(an-a+1)(p-1)\over ap-a+1}]}\frac{1}{\pi^{n-j}}
p^{[{(an-a+1)(p-1)-j(ap-a+1)-1 \over p-1}]}A^+.$$
\end{theorem}
\begin{proof} The theorem is true for $a=1$. Assume now that
$a>1$ and assume that the theorem is true for $a-1$. Then
\begin{eqnarray*}
&&\psi^a\left(\frac{1}{\pi^n}A^+\right)= \psi
\left( \psi^{a-1}\left(\frac{1}{\pi^n}A^+\right)\right)\\
& \subseteq & \psi \left( \sum_{j=0}^{[{(a-1)n-a+2)(p-1)\over
(a-1)p-a+2}]}\frac{1}{\pi^{n-j}}
p^{[{((a-1)n-a+2)(p-1)-j((a-1)p-a+2)-1\over p-1}]}A^+ \right)\\
& \subseteq & \sum_j \sum_i\frac{1}{\pi^{n-j-i}}
p^{[{((a-1)n-a+2)(p-1)-j((a-1)p-a+2)-1\over p-1}]+
[{(n-j)(p-1)-ip-1 \over p-1}]}A^+,
\end{eqnarray*}
where the indices $i$ and $j$ satisfy
$$0\leq j\leq [{(a-1)n-a+2)(p-1)\over
(a-1)p-a+2}], \ 0\leq i \leq [(n-j)(p-1)/p].
$$
For fixed $i+j=k$, the exponent of $p$ is decreasing in $j$ and
the minimun value is attained when $j=k$ and $i=0$. It follows
that
\begin{eqnarray*}
\psi^a\left(\frac{1}{\pi^n}A^+\right) &\subseteq&\sum_{k\geq
0}\frac{1}{\pi^{n-k}}
p^{[{((a-1)n-a+2)(p-1)-k((a-1)p-a+2)-1\over p-1}]+ [n-k-1]}A^+\\
&\subseteq&\sum_{k=0}^{[{(an-a+1)(p-1)\over ap-a+1}
]}\frac{1}{\pi^{n-k}} p^{[{(an-a+1)(p-1)-k(ap-a+1)-1\over
p-1}]}A^+,
\end{eqnarray*}
where we stop at $k=[{(an-a+1)(p-1)\over ap-a+1}]$ in the
summation as the exponent of $p$ is zero if $k\geq
[{(an-a+1)(p-1)\over ap-a+1} ]$.
\end{proof}

\begin{corollary}[Sharp Congruence II] Let
$$S_j(n,r,p^a)=\sum_{\begin{subarray}{c}i_0+\cdots+i_{p^a-1}=n\\
i_1+2i_2+\cdots \equiv r({\rm
mod}p^a)\end{subarray}}\binom{n}{i_0\cdots
i_{p^a-1}}\binom{(i_1+2i_2+\cdots -r)/p^a}{j}.$$ Then for integer
$j\geq 0$, we have
$$S_j(n,r,p^a)\equiv 0({\rm mod} p^{[{(an-a+1)(p-1)-j(ap-a+1)-1\over p-1}]}).$$
\end{corollary}

{}

\end{document}